\documentclass{article}%
\usepackage{latexsym}
\usepackage{amssymb}
\usepackage{amsfonts}
\usepackage{mathrsfs}
\usepackage{chemarrow}
\usepackage{amsthm,array,bm}
\usepackage{pgf}
\usepackage{pgflibrarysnakes}
\usepackage{tikz}
\usepackage{mathtools}
\usepackage{array}
\usepackage{color,graphicx}
\usepackage{pstricks}
\usepackage{amsmath}
\usepackage{graphicx}%
\providecommand{\U}[1]{\protect\rule{.1in}{.1in}}
\def\bz{\mathbb Z}

\newtheorem{theorem}{Theorem}[section]

\newtheorem{lemma}[theorem]{Lemma}

{\bfseries}{\rmfamily}
\newtheorem{definition}[theorem]{Definition}

\newtheorem{remark-definition}[theorem]{Remark and Definition}
\font\smc=cmcsc10
\newrgbcolor{lightblueA}{.50 1 1}
\newrgbcolor{violet}{.6 .1 .8} 
\newrgbcolor{lightyellow}{1 1 .8} 
\newrgbcolor{lightblue}{.80 1 1}
\newrgbcolor{mygreen}{0 .66 .05} 
\definecolor{mygreen}{rgb}{0,.66,.05}
\definecolor{lightyellow}{rgb}{1,1,.80}
\newrgbcolor{orange}{1 .6 0}
\newrgbcolor{GreenYellow}{.85 1 .31}
\newrgbcolor{Yellow}{1  1  0}
\newrgbcolor{Goldenrod}{1  .90  .16}
\newrgbcolor{Dandelion}{1  .71  .16}
\newrgbcolor{Apricot}{1  .68  .48}
\newrgbcolor{Peach}{1  .50  .30}
\newrgbcolor{Melon}{1  .54  .50}
\newrgbcolor{YellowOrange}{1  .58  0}
\newrgbcolor{Orange}{1  .39  .13}
\newrgbcolor{BurntOrange}{1  .49  0}
\newrgbcolor{Bittersweet}{1.  .4300  .24}
\newrgbcolor{RedOrange}{1  .23  .13}
\newrgbcolor{Mahogany}{1.  .4475  .4345}
\newrgbcolor{Maroon}{1.  .4084  .5376}
\newrgbcolor{BrickRed}{1.  .3592  .3232}
\newrgbcolor{Red}{1  0  0}
\newrgbcolor{OrangeRed}{1  0  .50}
\newrgbcolor{RubineRed}{1  0  .87}
\newrgbcolor{WildStrawberry}{1  .04  .61}
\newrgbcolor{Salmon}{1  .47  .62}
\newrgbcolor{CarnationPink}{1  .37  1}
\newrgbcolor{Magenta}{1  0  1}
\newrgbcolor{VioletRed}{1  .19  1}
\newrgbcolor{Rhodamine}{1  .18  1}
\newrgbcolor{Mulberry}{.6668  .1180  1.}
\newrgbcolor{RedViolet}{.9538  .4060  1.}
\newrgbcolor{Fuchsia}{.5676  .1628  1.}
\newrgbcolor{Lavender}{1  .52  1}
\newrgbcolor{Thistle}{.88  .41  1}
\newrgbcolor{Orchid}{.68  .36  1}
\newrgbcolor{DarkOrchid}{.60  .20  .80}
\newrgbcolor{Purple}{.55  .14  1}
\newrgbcolor{Plum}{.50  0  1}
\newrgbcolor{Violet}{.98 .15 .95}
\newrgbcolor{RoyalPurple}{.25  .10  1}
\newrgbcolor{BlueViolet}{.84  .38  .98}
\newrgbcolor{Periwinkle}{.43  .45  1}
\newrgbcolor{CadetBlue}{.38  .43  .77}
\newrgbcolor{CornflowerBlue}{.35  .87  1}
\newrgbcolor{MidnightBlue}{.4414  .9259  1.}
\newrgbcolor{NavyBlue}{.06  .46  1}
\newrgbcolor{RoyalBlue}{0  .50  1}
\newrgbcolor{Blue}{0  0  1}
\newrgbcolor{Cerulean}{.06  .89  1}
\newrgbcolor{Cyan}{0  1  1}
\newrgbcolor{ProcessBlue}{.04  1  1}
\newrgbcolor{SkyBlue}{.38  1  .88}
\newrgbcolor{Turquoise}{.15  1  .80}
\newrgbcolor{TealBlue}{.1572  1.  .6668}
\newrgbcolor{Aquamarine}{.18  1  .70}
\newrgbcolor{BlueGreen}{.15  1  .67}
\newrgbcolor{Emerald}{0  1  .50}
\newrgbcolor{JungleGreen}{.01  1  .48}
\newrgbcolor{SeaGreen}{.31  1  .50}
\newrgbcolor{Green}{0  1  0}
\newrgbcolor{ForestGreen}{.1992  1.  .2256}
\newrgbcolor{PineGreen}{.3100  1.  .5575}
\newrgbcolor{LimeGreen}{.50  1  0}
\newrgbcolor{YellowGreen}{.56  1  .26}
\newrgbcolor{SpringGreen}{.74  1  .24}
\newrgbcolor{OliveGreen}{.6160  1.  .4300}
\newrgbcolor{RawSienna}{.53  .28  .16}
\newrgbcolor{Sepia}{1.  .7510  .70}
\newrgbcolor{Brown}{.41  .25  .18}
\newrgbcolor{TAN}{.86  .58  .44}
\newrgbcolor{Gray}{1.  1.  1.}
\newrgbcolor{Black}{1  1  1}
\newrgbcolor{White}{1  1  1}
\begin{document}

\title{Symmetries of Nonlinear Vibrations in Tetrahedral Molecular Configurations}
\author{Irina Berezovik\thanks{Department of Mathematical Sciences University of Texas
at Dallas Richardson, 75080 USA} , Carlos Garc\'{\i}a-Azpeitia\thanks{
Departamento de Matem\'{a}ticas, Facultad de Ciencias, Universidad Nacional
Aut\'{o}noma de M\'{e}xico, 04510 M\'{e}xico DF, M\'{e}xico.
cgazpe@ciencias.unam.mx},\\and Wieslaw Krawcewicz$^{*}$\thanks{Center for Applied Mathematics, Guangzhou
University, Guangzhou, China} }
\maketitle

\begin{abstract}
We study nonlinear vibrational modes of oscillations for tetrahedral
configurations\ of particles. In the case of tetraphosphorus, the interaction
of atoms is given by bond stretching and van der Waals forces. Using
equivariant gradient degree, we present a topological classification of the
spatio-temporal symmetries of the periodic solutions with finite Weyl's group.
This procedure gives all the symmetries of the nonlinear vibrations for
general force fields.

MSC 37J45, 34C25, 37G40, 47H11, 70H33

\end{abstract}

\section{Introduction}

Description of the dynamical motions of a collection of particles in space and
time can provide a rich amount of information including molecular geometries,
mean atomic fluctuations, and free energies. The molecular conformation is
located at a local energy minimum where the net inter-particle force on each
particle is close to zero and the position on the potential energy surface is
stationary. Molecular motion can be modeled as vibrations around and
interconversions between these stable configurations. Molecules spend most of
their time in these low-lying states at finite temperature, which thus
dominate the molecular properties of the system.

In this paper we study the molecular mechanics of tetrahedral molecules. Let
$u(t)=(u_{1}(t),\,u_{2}(t),\,u_{3}(t),\,u_{4}(t))$ with $u_{j}(t)\in
\mathbb{R}^{3}$ for $j=1,2,3,4$ stand for the spatial position of the system
of $4$-particles at time $t$. Such system satisfies the Newtonian equation
\begin{equation}
\ddot{u}(t)=-\nabla V(u(t)), \label{eqn01}%
\end{equation}
where the potential energy $V$ represents the force field given by
\[
V(u):=\sum_{1\leq j<k\leq4}^{n}U(|u_{j}-u_{k}|^{2}).
\]
When these $4$-particles interact by bond stretching, van der Waals and
electrostatic forces \cite{xx, GT}, $U$ is given by
\[
U(x)=\left(  \sqrt{x}-1\right)  ^{2}+\left(  \frac{B}{x^{6}}-\frac{A}{x^{3}%
}\right)  +\frac{\sigma}{\sqrt{x}}~.
\]

A local energy minimum is a stationary point $a\in\mathbb{R}^{12}$ such that
$\nabla V(a)=0$. To detect possible periodic vibrations around the
configuration $a$, a natural method is to investigate the existence of
periodic solutions to \eqref{eqn01} near $a$. An important feature of this
molecular configuration is that it admits tetrahedral spatial symmetries and
thus the bifurcated/emerging periodic motions will have both spatial and
temporal symmetry. In this case, the equation \eqref{eqn01} is equivariant
under the action of the group
\[
S_{4}\times O(3)\times O(2)\text{,}%
\]
which acts by permuting the particles, rotating and reflecting them in
$\mathbb{R}^{3}$ and by temporal phase shift and reflection, respectively.

In this paper, we use the equivariant degree method to investigate the
existence of periodic solution to \eqref{eqn01} around an equilibrium
admitting $S_{4}$-symmetries. The concept of equivariant gradient degree was
introduced by K. Geba in \cite{Geba}. This degree satisfies all the standard
properties expected from a degree theory. In addition, it can also be
generalized to settings in infinite-dimensional spaces allowing its
applications to studying critical points of invariant functionals (cf.
\cite{survey}). The values of the gradient equivariant degree can be expressed
elegantly in the form
\[
\nabla_{G}\text{-deg}(\nabla\varphi,\Omega)=n_{1}(H_{1})+n_{2}(H_{2}%
)+\dots+n_{m}(H_{m}),\;\;\;n_{k}\in\mathbb{Z},
\]
where $(H_{j})$ are the orbit types in $\Omega$, which allow to predict the
existence of various critical orbits for $\varphi$ and their symmetries. We
should mention that the gradient degree is just one of many equivariant
degrees (see \cite{survey}) that were introduced in the last three decades:
equivariant degrees with $n$-free parameters (primary degrees, twisted
degrees), gradient and orthogonal equivariant degrees \cite{AED, IVB, KW} ---
all these different degrees being related to each other (cf. \cite{BKR,RR}).
For multiple applications of the equivariant gradient degree to Newtonian
system, we refer the reader to \cite{survey, DKY, FRR, GI, RR, RY2} and the
references therein.

The local minimizer $u_{o}$ of $V$ is a regular tetrahedron located in a
sphere of radius $r_{o}$. Then $U^{\prime\prime}(r_{o})>0$ due to the fact
that $a$ is a local minimizer. Let%
\[
\nu_{0}^{2}:=\frac{32}{3}r_{o}^{2}U^{\prime\prime}(r_{o})>0\text{.}%
\]
The 6 non-zero eigenvalues of $D^{2}V(u_{o})$ are computed to be $\nu_{0}^{2}$
with multiplicity 2, $2\nu_{0}^{2}$ with multiplicity 3{,} and $4\nu_{0}^{2}$
with multiplicity 1. Then, the normal modes of (\ref{eqn01}) are $\nu_{0}$,
$\sqrt{2}\nu_{0}$, and $2\nu_{0}$ with the respective multiplicities.

Observe that the normal mode $\nu_{0}$ is $1:1:2$ resonant. Due to
multiplicities and resonances, the Lyapunov center theorem can be applied to
prove only the local existence of a periodic solution (nonlinear normal modes)
from the frequency $2\nu_{0}$ \cite{MW}. On the other hand, since the
equilibrium corresponds to a local minimizer of the Hamiltonian, the
Weinstein-Moser theorem \cite{1} gives the existence of at least 6 periodic
orbits in each (small) fixed energy level, regardless of resonances and
multiplicities. Using the gradient equivariant degree method, we establish the
global existence of branches of periodic solutions emerging from the
equilibrium $u_{o}$ starting with the frequencies of the normal modes $\nu
_{0}$, $\sqrt{2}\nu_{0}$, and $2\nu_{0}$. The global property means that
families of periodic solutions are represented by a continuum that has norm or
period going to infinity, ends in a collision orbit, or comes back to another equilibrium.

Specifically, we prove that the tetrahedral equilibrium $u_{o}$ has the
following global family of periodic solutions: one family starting with
frequency $\nu_{0}$, five families with frequency $\sqrt{2}\nu_{0}$, and one
families with frequency $2\nu_{0}$. The family from $\nu_{0}$ has symmetries
of a brake orbits where all the particles form a regular tetrahedron at any
time. The first symmetry from $\sqrt{2}\nu_{0}$ gives brake orbits where two
pairs of particles are related by inversion, and in the second symmetry, one
pair of particles is related by inversion and other by a $\pi$-rotation and
$\pi$-phase shift. The third symmetry from $\sqrt{2}\nu_{0}$ is not a brake
orbit, and all the particles are related by a $\pi/2$-rotoreflection and
$\pi/2$-phase shift. The fourth symmetry from $\sqrt{2}\nu_{0}$ is a brake
orbit where three particles form a triangle at all times, while another makes
counterbalance movement. The fifth symmetry from $\sqrt{2}\nu_{0}$ is not
brake orbit, and three particles move in the form of a traveling wave along a
triangle, while another makes a counterbalance movement. The family from
$2\nu_{0}$ has symmetries of a brake orbits with two symmetries by inversion
at any time. The exact description of the symmetries is given in Section 5.

The article \cite{5} presents for tetraphosphorus molecules an extensive study
of the stability and existence of nonlinear modes that are relative
equilibria. The authors assume the absence of resonances in the normal form of
that Hamiltonian. Though in our study we consider the nonlinear normal modes
that are not relative equilibria and we use a force field expressing mutual
interaction between the atoms that leads to a Hamiltonian with resonances.

Other molecules that have tetrahedral symmetries include the methane molecule.
This molecule has an equilibrium state with a carbon atom at the center and
four hydrogen atoms at the vertices of a regular tetrahedron. The articles
\cite{3,4} use a combination of geometric methods, normal forms, and Krein
signature to analyze the existence of nonlinear modes and their stability.
These results can be easily extrapolated to the tetraphosphorus molecule which
have the same symmetries but different configuration. In this sense, the
symmetries and number of solutions obtained in \cite{3,4} for each frequency
coincide with our results. The gradient equivariant degree allows to determine
global properties of the branches and to manage easily resonances.
Nevertheless, more precise local information can be obtain with the results of
\cite{3,4,5}.

The paper is structured in the following sections. In Section 2, we analyze
the isotypic decomposition of the eigenvalues of the Hessian $D^{2}V(u_{o})$.
In Section 3, we prove the global existence of families of periodic solutions
from the tetrahedral equilibrium. In Section 4, we describe the symmetries of
the different families of periodic solutions. In Appendix, we review
preliminary notions and definitions used in group representations, the
properties of the equivariant gradient degree, and indicate the standard
techniques used to compute it.

\section{Model for Atomic Interaction}

Consider $4$ identical particles $u_{j}$ in the space $\mathbb{R}^{3}$, for
$j=1,2,3,4$. Assume that each particle $u_{j}$ interacts with all other
particles $u_{k}$ for $k\not =j$. Put $u:=(u_{1},u_{2},u_{3},u_{4})^{T}%
\in\mathbb{R}^{12}$ and
\[
\tilde{\Omega}_{o}:=\{u\in\mathbb{R}^{12}:\forall_{k\not =j}\;\;u_{k}%
\not =u_{j}\}.
\]
The Newtonian equation that describes the interaction between these
$4$-particles is%
\begin{equation}
\ddot{u}=-\nabla V(u),\quad u\in\tilde{\Omega}_{o}. \label{eq:mol}%
\end{equation}

The potential energy $V:\tilde{\Omega}_{o}\rightarrow\mathbb{R}$,%
\begin{equation}
V(u):=\sum_{1\leq j<k\leq4}^{n}U(|u_{j}-u_{k}|^{2}), \label{eq:pot}%
\end{equation}
is well define and $C^{2}$ when $U\in C^{2}(\mathbb{R}^{+})$ satisfies%
\begin{equation}
\lim_{x^{+}\rightarrow0}U(x)=\infty,\qquad\lim_{x\rightarrow\infty}%
U(x)=\infty\text{.} \label{U}%
\end{equation}
Classical forces used in molecular mechanics are associated with bending
between adjacent particles, electrostatic interactions and van der Walls
forces. The condition (\ref{U}) holds when $U$ is determined by these force fields.

\subsection{The Tetrahedral Equilibrium}

\label{sec:equilib} One can easily notice that the space $\mathbb{R}^{12}$ is
a representation of the group
\[
\mathfrak{G}:=S_{4}\times O(3),
\]
where $S_{4}$ stands for the symmetric group of four elements. More precisely
$S_{4}$ is the group of permutations of four elements $\{1,2,3,4\}$. Then the
action of $\mathfrak{G}$ on $\mathbb{R}^{12}$ is given by
\begin{equation}
(\sigma,A)(u_{1},u_{2},u_{3},u_{4})^{T}=(Au_{\sigma(1)},Au_{\sigma
(2)},Au_{\sigma(3)},Au_{\sigma(4)})^{T}, \label{eq:act1}%
\end{equation}
where $A\in O(3)$ and $\sigma\in S_{4}$.

Notice that $S_{4}$ can be considered as a subgroup of $O(3)$, representing
the actual symmetries of a tetrahedron $\mathbf{T}\subset\mathbb{R}^{3}$. More
precisely, consider the regular tetrahedron given by
\[
\mathbf{T}:=\left\{  \gamma_{1},\gamma_{2},\gamma_{3},\gamma_{4}\right\}  ,
\]
where
\[
\gamma_{1}=\left(
\begin{array}
[c]{c}%
0\\
0\\
1
\end{array}
\right)  ,\quad\gamma_{2}=\left(
\begin{array}
[c]{c}%
\frac{2}{3}\sqrt{2}\\
0\\
-\frac{1}{3}%
\end{array}
\right)  ,\quad\gamma_{3}=\left(
\begin{array}
[c]{c}%
-\frac{1}{3}\sqrt{2}\\
\frac{1}{3}\sqrt{6}\\
-\frac{1}{3}%
\end{array}
\right)  ,\quad\gamma_{4}=\left(
\begin{array}
[c]{c}%
-\frac{1}{3}\sqrt{2}\\
-\frac{1}{3}\sqrt{6}\\
-\frac{1}{3}%
\end{array}
\right)  .
\]
The tetrahedral group $\{A\in O(3):A(\mathbf{T})=\mathbf{T}\}$ can be
identified with the group $S_{4}$. Indeed, any $A$ such that $A(\mathbf{T}%
)=\mathbf{T}$ permutes the vertices of $\mathbf{T}$, i.e.
\[
A\gamma_{j}=\gamma_{\sigma(j)}%
\]
for $j=1,2,3,4$, i.e. we can identify $A_{\sigma}$ with the permutation
$\sigma\in S_{4}$ by these relations. We have explicitly for the permutations
$(1,2)$ and $(2,3,4)$, which are generators of $S_{4}$, the following
identification
\[
A_{(1,2)}=\left[
\begin{array}
[c]{ccc}%
\frac{1}{3} & 0 & \frac{2\sqrt{2}}{3}\\
0 & 1 & 0\\
\frac{2\sqrt{2}}{3} & 0 & -\frac{1}{3}%
\end{array}
\right]  \quad\text{ and }\quad A_{(2,3,4)}=\left[
\begin{array}
[c]{ccc}%
-\frac{1}{2} & \frac{\sqrt{3}}{2} & 0\\
-\frac{\sqrt{3}}{2} & -\frac{1}{2} & 0\\
0 & 0 & 1
\end{array}
\right]  ~.
\]
These generators define an explicit isomorphism $A_{\sigma}:S_{4}\rightarrow
O(3)$.

Notice that the function $V:\tilde{\Omega}_{o}\rightarrow\mathbb{R}$ is
invariant with respect to the action of $c\in\mathbb{R}^{3}$ on $(\mathbb{R}%
^{3})^{4}$ by shifting, $V(u+c)=V(u)$. Therefore in order to fix the center of
mass at the origin in the system (\ref{eq:mol}), we define the subspace
\begin{equation}
\mathscr V:=\{(u_{1},u_{2},u_{3},u_{4})^{T}\in(\mathbb{R}^{3})^{4}:u_{1}%
+u_{2}+u_{3}+u_{4}=0\} \label{eq:V}%
\end{equation}
and $\Omega_{o}=\tilde{\Omega}_{o}\cap\mathscr V$. Then, one can easily notice
that $\mathscr V$ and $\Omega_{o}$ are\emph{ }invariant under the nonlinear
dynamics of \eqref{eq:mol}, and in addition $\Omega_{o}$ is $G$-invariant.

\begin{figure}[ptb]
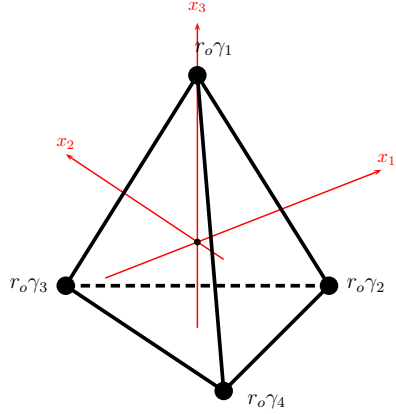

\vglue2.1cm\hskip3cm \scalebox{.7}{
\psline[linecolor=red]{->}(2.5,-.8)(2.5,5)
\rput(2.5,5.3){\red $x_3$}
\psline[linecolor=red]{->}(0.75,0.145)(6,2.2)
\rput(6.1,2.4){\red $x_1$}
\rput(0,2.5){\psline[linecolor=red]{<-}(0,0)(3,-2)}
\rput(0,2.7){\red $x_2$}
\psdot(2.5,.83)
\psline[linecolor=white,fillcolor=lightblue](0,0)(3,-2)(2.5,4)
\psline[linecolor=white,fillcolor=lightyellow](2.5,4)(3,-2)(5,0)
\psline[linewidth=2pt](0,0)(3,-2)(5,0)
\psline[linewidth=2pt](0,0)(2.5,4)(5,0)
\psline[linewidth=2pt](3,-2)(2.5,4)
\psline[linewidth=2pt,linestyle=dashed](0,0)(5,0)
\psdot[dotsize=10pt](5,0)\psdot[dotsize=10pt](0,0)\psdot[dotsize=10pt](3,-2)\psdot[dotsize=10pt](2.5,4)
\rput(5.7,0){\large $r_{o}\gamma_2$}
\rput(-.7,0){\large $r_{o}\gamma_3$}\
\rput(3.7,-2.2){\large $r_{o}\gamma_4$}
\rput(2.7,4.5){\large $r_{o}\gamma_1$}}
\par
\vskip1.5cm \caption{Stationary solution to equation \eqref{eq:mol} with
tetrahedral symmetries.}%
\label{fig:1}%
\end{figure}

Consider the point $v_{0}:=(\gamma_{1},\gamma_{2},\gamma_{3},\gamma_{4}%
)\in\Omega_{o}$. The isotropy group $\mathfrak{G}_{v_{o}}$ is given by
\[
\tilde{S}_{4}:=\{(\sigma,A_{\sigma})\in S_{4}\times O(3):\sigma\in S_{4}\},
\]
where $S_{4}$ is considered as a subgroup of $O(3)$ using the above
identification for $A_{\sigma}$. Since $\tilde{S}_{4}$ is a finite group,
$\mathscr V^{\tilde{S}_{4}}$ is a one dimensional subspace of $\mathscr V$ and
we have that
\[
\mathscr V^{\tilde{S}_{4}}=\text{span}_{\mathbb{R}}\{(\gamma_{1},\gamma
_{2},\gamma_{3},\gamma_{4})^{T}\}.
\]
Then, by Symmetric Criticality Condition, a critical point of $V^{\tilde
{S}_{4}}:\Omega_{o}^{\tilde{S}_{4}}\rightarrow\mathbb{R}$ is also a critical
point of $V$. Since $\mathscr V^{\tilde{S}_{4}}$ is one-dimensional, we denote
its vectors by $rv_{o}\in\mathbb{R}^{12}$ for $r\in\mathbb{R}$. Notice that
\[
\phi(r):=\sum_{1\leq j<k\leq4}U\left(  \frac{8}{3}r^{2}\right)  ,\quad r>0.
\]
is exactly the restriction of $V$ to the fixed-point subspace $\mathscr
V^{\tilde{S}_{4}}\cap\Omega_{o}$. Thus in order to find an equilibrium for
\eqref{eq:mol}, by Symmetric Criticality Principle, it is sufficient to
identify a critical point $r_{o}$ of $\phi(r)$. Clearly by (\ref{U}),
\[
\lim_{r\rightarrow0^{+}}\phi(r)=\lim_{r\rightarrow\infty}\phi(r)=\infty,
\]
then there exists a minimizer $r_{o}\in(0,\infty)$, which is clearly a
critical point of $\varphi$. Consequently
\begin{equation}
u_{o}:=r_{o}v_{o}\in\Omega_{o} \label{eq:crit-u0}%
\end{equation}
is the $\tilde{S}_{4}$-symmetric equilibrium of $V$. The components of $u_{o}%
$, which are $r_{o}\gamma_{j}$ for $j=1,2,3,4$, give us the configuration of
the stationary solution of \eqref{eq:mol}, see Figure \ref{fig:1}.

\subsection{Isotypic Decomposition}

Since the system \eqref{eq:mol} is symmetric with respect to the group action
$\mathfrak{G}:=S_{4}\times O(3)$ we have that the orbit of equilibria
$\mathfrak{G}(u_{o})$ is a 3-dimensional submanifold in $\mathscr V$. The
slice $S_{o}$ to the orbit $\mathfrak{G}(u_{o})$ at $u_{o}$ is
\[
S_{o}:=\{x\in\mathscr V:x\bullet T_{u_{o}}\mathfrak{G}(u_{o})=0\}.
\]
The tangent space $T_{u_{o}}\mathfrak{G}(u_{o})$ is described as%
\[
T_{u_{o}}\mathfrak{G}(u_{o})=\text{span}\{(J_{j}\gamma_{1},J_{j}\gamma
_{2},J_{j}\gamma_{3},J_{j}\gamma_{4})^{T}\in\mathscr V:j=1,2,3\},
\]
where $J_{j}$ are the three infinitesimal generator of the rotations:%
\[
J_{1}:=\left[
\begin{array}
[c]{ccc}%
0 & 0 & 0\\
0 & 0 & -1\\
0 & 1 & 0
\end{array}
\right]  ,\quad J_{2}:=\left[
\begin{array}
[c]{ccc}%
0 & 0 & 1\\
0 & 0 & 0\\
-1 & 0 & 0
\end{array}
\right]  ,\quad J_{3}:=\left[
\begin{array}
[c]{ccc}%
0 & -1 & 0\\
1 & 0 & 0\\
0 & 0 & 0
\end{array}
\right]  .
\]

Since the $\mathfrak{G}$-isotropy group of $u_{o}$ is $\tilde{S}_{4}$, then
$S_{o}$ is an orthogonal $\tilde{S}_{4}$ representation. In order to identify
the $\tilde{S}_{4}$-isotypical components, we consider first the $\tilde
{S}_{4}$-representation $V=\mathbb{R}^{12}$ on which $\tilde{S}_{4}$-acts by
\eqref{eq:act1}. We have the following table of characters $\chi_{j}$,
$j=0,1,2,3,4$, for all irreducible $\tilde{S}_{4}$-representations
$\mathcal{V}_{j}$ (all of them of real type) and the character $\chi_{V}$ of
the the representation $V$:%

\[%
\begin{tabular}
[c]{|cc|ccccc|}\hline
Rep. & Character & $(1)$ & $(1,2)$ & $(1,2)(3,4)$ & $(1,2,3)$ & $(1,2,3,4)$%
\\\hline
$\mathcal{V}_{0}$ & $\chi_{0}$ & $1$ & $1$ & $1$ & $1$ & $1$\\
$\mathcal{V}_{1}$ & $\chi_{1}$ & $3$ & $1$ & $-1$ & $0$ & $-1$\\
$\mathcal{V}_{2}$ & $\chi_{2}$ & $2$ & $0$ & $2$ & $-1$ & $0$\\
$\mathcal{V}_{3}$ & $\chi_{3}$ & $3$ & $-1$ & $-1$ & $0$ & $1$\\
$\mathcal{V}_{4}$ & $\chi_{4}$ & $1$ & $-1$ & $1$ & $1$ & $-1$\\\hline
$V$ & $\chi_{V}$ & $12$ & $2$ & $0$ & $0$ & $0$\\\hline
\end{tabular}
\]

One can easily conclude that we have the following $\tilde{S}_{4}$-isotypic
decomposition:
\[
V=\mathcal{V}_{0}\oplus\left(  \mathcal{V}_{1}\oplus\mathcal{V}_{1}\right)
\oplus\mathcal{V}_{2}\oplus\mathcal{V}_{3}~.
\]
Since the subspace $V$ is obtained by fixing the center of mass at the origin,
and $\{(v,v,v,v)\in\mathbb{R}^{12}:v\in\mathbb{R}^{3}\}$ is equivalent to the
irreducible $\tilde{S}_{4}$-representation $\mathcal{V}_{1}$, we have
the\emph{ }the $\tilde{S}_{4}$-isotypic decomposition%
\begin{equation}
\mathscr V=V_{0}\oplus V_{1}\oplus V_{2}\oplus V_{3},\qquad V_{j}%
=\mathcal{V}_{j}~. \label{eq:S4-iso}%
\end{equation}

In order to determine the $\tilde{S}_{4}$-isotypic type of the tangent space
$T_{u_{0}}\mathfrak{G}(u_{o})$ (which has to be an irreducible $\tilde{S}_{4}%
$-representation of dimension $3$), we apply the isotypic projections
$P_{j}:V\rightarrow V_{j}$, $j=1$ and $3$, given by
\[
P_{j}v:=\frac{\dim(V_{j})}{72}\sum_{g\in S_{4}}\chi_{j}(g)\,gv,\quad v\in V,
\]
to conclude that $T_{u_{0}}\mathfrak{G}(u_{o})\simeq\mathcal{V}_{3}$.
Therefore, the $\tilde{S}_{4}$-isotypic decomposition of the slice $S_{o}$ is
\begin{equation}
S_{o}=V_{0}\oplus V_{1}\oplus V_{2}~.\label{eq:S4-iso-S}%
\end{equation}

\subsection{Computation of the Spectrum $\sigma(\nabla^{2}V(u_{o}))$}

Since the potential $V$ is given by \eqref{eq:pot}, we have%

\[
\nabla V(u)=2\left[
\begin{array}
[c]{c}%
\sum_{k\not =0}U^{\prime}(|u_{0}-u_{k}|^{2})(u_{0}-u_{k})\\
\sum_{k\not =1}U^{\prime}(|u_{1}-u_{k}|^{2})(u_{1}-u_{k})\\
\vdots\\
\sum_{k\not =n-1}U^{\prime}(|u_{n-1}-u_{k}|^{2})(u_{n-1}-u_{k}).
\end{array}
\right]
\]
Notice that we have $\nabla V(u_{o})=0$ when $U^{\prime}(r_{0}^{2})=0$.

For a given vector $v=(x,y,z)^{T}\in\mathbb{R}^{3}$, we define the matrix
$\mathfrak{m}_{v}:=vv^{T}$, i.e.
\[
\mathfrak{m}_{v}:=\left[
\begin{array}
[c]{c}%
x\\
y\\
z
\end{array}
\right]  [x,y,z]=\left[
\begin{array}
[c]{ccc}%
x^{2} & xy & xz\\
xy & y^{2} & yz\\
xz & yz & z^{2}%
\end{array}
\right]  .
\]
Then one can easily see that the matrix $\mathfrak{m}_{v}$ represents the
linear operator $\Vert v\Vert^{2}\,P_{v}:\mathbb{R}^{3}\rightarrow
\mathbb{R}^{3}$, where $P_{v}$ is the orthogonal projection onto the subspace
generated by $v\in\mathbb{R}^{3}$. Put
\[
\mathfrak{m}_{j,k}:=\mathfrak{m}_{(\gamma_{j}-\gamma_{k})}.
\]
Clearly $\mathfrak{m}_{j,k}=\mathfrak{m}_{k,j}$. Notice that
\[
\mathfrak{m}_{j,k}(\gamma_{j})=\frac{4}{3}(\gamma_{j}-\gamma_{k}).
\]

By direct computations one can derive the following matrix form of $\nabla
^{2}V(u_{o})$
\[
M:=\nabla^{2}V(u_{o})=4r_{o}^{2}U^{\prime\prime}(r_{o})\left[
\begin{array}
[c]{cccc}%
\displaystyle\sum_{j\not =1}\mathfrak{m}_{1j} & -\mathfrak{m}_{12} &
-\mathfrak{m}_{13} & -\mathfrak{m}_{1,4}\\
-\mathfrak{m}_{2,1} & \displaystyle\sum_{j\not =2}\mathfrak{m}_{2,j} &
-\mathfrak{m}_{23} & -\mathfrak{m}_{2,4}\\
-\mathfrak{m}_{3,1} & -\mathfrak{m}_{3,2} & \displaystyle\sum_{j\not =%
3}\mathfrak{m}_{3,j} & -\mathfrak{m}_{3,4}\\
-\mathfrak{m}_{4,1} & -\mathfrak{m}_{3,2} & -\mathfrak{m}_{4,3} &
\displaystyle\sum_{j\not =4}\mathfrak{m}_{4,j}%
\end{array}
\right]
\]
Since $M:\mathscr V\rightarrow\mathscr V$ is $\tilde{S}_{4}$-equivariant, it
follows that
\[
M_{j}:=M|_{V_{j}}:V_{j}\rightarrow V_{j},\quad j=0,1,2,
\]
and since the sub-representations $V_{j}=\mathcal{V}_{j}$ are absolutely
irreducible, we have that
\[
M_{j}=\mu_{j}\,\text{\textrm{Id\thinspace}}:V_{j}\rightarrow V_{j}%
,\;\;\;j=0,1,2,
\]
which implies $\sigma(M|_{{S_{o}}})=\{\mu_{0},\mu_{1},\mu_{2})$.

In order to find explicit formulae for the eigenvalues $\mu_{j}$, we notice
that
\[
\mathfrak{v}_{0}:=\left[
\begin{array}
[c]{c}%
\gamma_{1}\\
\gamma_{2}\\
\gamma_{3}\\
\gamma_{4}%
\end{array}
\right]  \in V_{0},\;\;\mathfrak{v}_{1}:=\left[
\begin{array}
[c]{c}%
-2\gamma_{1}\\
\gamma_{1}+\gamma_{2}\\
\gamma_{1}+\gamma_{3}\\
\gamma_{1}+\gamma_{4}%
\end{array}
\right]  \in V_{1},\;\;\mathfrak{v}_{2}:=\left[
\begin{array}
[c]{c}%
\gamma_{2}-\gamma_{3}\\
\gamma_{1}-\gamma_{4}\\
\gamma_{4}-\gamma_{1}\\
\gamma_{3}-\gamma_{2}%
\end{array}
\right]  \in V_{2},
\]
and by direct application of the matrix $\mathscr L$ on the vectors
$\mathfrak{v}_{j}$, $j=0,1,2$, we obtain that
\begin{align*}
\mu_{0}  &  =\frac{128}{3}r_{o}^{2}U^{\prime\prime}(r_{o})=4\nu_{0}^{2},\\
\mu_{1}  &  =\frac{64}{3}r_{o}^{2}U^{\prime\prime}(r_{o})=2\nu_{0}^{2},\\
\mu_{2}  &  =\frac{32}{3}r_{o}^{2}U^{\prime\prime}(r_{o})=\nu_{0}^{2}.
\end{align*}
Notice that $0<\mu_{2}<\mu_{1}<\mu_{0}~$.

\section{Equivariant Bifurcation}

In what follows, we are interested in finding non-trivial $T$-periodic
solutions to \eqref{eq:mol}, bifurcating from the orbit of equilibrium points
$\mathfrak{G}(u_{o})$. By normalizing the period, i.e. by making the
substitution $v(t):=u\left(  \frac{T}{2\pi}t\right)  $ in \eqref{eq:mol}, we
obtain the following system
\begin{equation}%
\begin{cases}
\ddot{v}=-\lambda^{2}\nabla V(v),\\
v(0)=v(2\pi),\;\;\dot{v}(0)=\dot{v}(2\pi),
\end{cases}
\label{eq:mol1}%
\end{equation}
where $\lambda^{-1}=2\pi/T$ is the frequency.

\subsection{Equivariant Gradient Map}

Since $\mathscr V$ is an orthogonal $\mathfrak{G}$- representation, we can
consider the first Sobolev space of $2\pi$-periodic functions from
$\mathbb{R}$ to $\mathscr V$, i.e.
\[
H_{2\pi}^{1}(\mathbb{R},\mathscr V):=\{u:\mathbb{R}\rightarrow\mathscr
V\;:\;u(0)=u(2\pi),\;u|_{[0,2\pi]}\in H^{1}([0,2\pi];\mathscr V)\},
\]
equipped with the inner product
\[
\langle u,v\rangle:=\int_{0}^{2\pi}(\dot{u}(t)\bullet\dot{v}(t)+u(t)\bullet
v(t))dt~.
\]

Let $O(2)=SO(2)\cup\kappa SO(2)$ denote the group of $2\times2$-orthogonal
matrices, where
\[
\kappa=\left[
\begin{array}
[c]{cc}%
1 & 0\\
0 & -1
\end{array}
\right]  ,\qquad\left[
\begin{array}
[c]{cc}%
\cos\tau & -\sin\tau\\
\sin\tau & \cos\tau
\end{array}
\right]  \in SO(2)~.
\]
It is convenient to identify a rotation with $e^{i\tau}\in S^{1}%
\subset\mathbb{C}$. Notice that $\kappa e^{i\tau}=e^{-i\tau}\kappa$, i.e
$\kappa$ as a linear transformation of $\mathbb{C}$ into itself, acts as
complex conjugation.

Clearly, the space $H_{2\pi}^{1}(\mathbb{R},\mathscr V)$ is an orthogonal
Hilbert representation of
\[
G:=\mathfrak{G}\times O(2),\qquad\mathfrak{G}=S_{4}\times O(3).
\]
Indeed, we have for $u\in H_{2\pi}^{1}(\mathbb{R},\mathscr V)$ and
$(\sigma,A)\in\mathfrak{G}$ (see \eqref{eq:act1})
\begin{align}
\left(  \sigma,A\right)  u(t)  &  =(\sigma,A)u(t),\label{eq:ac}\\
e^{i\tau}u(t)  &  =u(t+\tau),\nonumber\\
\kappa u(t)  &  =u(-t).\nonumber
\end{align}

It is useful to identify a $2\pi$-periodic function $u:\mathbb{R}\rightarrow
V$ with a function $\widetilde{u}:S^{1}\rightarrow\mathscr V$ via the map
{$\mathfrak{e}(\tau)=e^{i\tau}:\mathbb{R}$}$\rightarrow S^{1}$. Using this
identification, we will write $H^{1}(S^{1},\mathscr V)$ instead of $H_{2\pi
}^{1}(\mathbb{R},\mathscr V)$.

Put
\[
\Omega:=\{u\in H^{1}(S^{1},\mathscr V):u(t)\in\Omega_{o}\}.
\]
Then, the system \eqref{eq:mol1} can be written as the following variational
equation
\begin{equation}
\nabla_{u}J(\lambda,u)=0,\quad(\lambda,u)\in\mathbb{R}\times\Omega,
\label{eq:bif1}%
\end{equation}
where $J:\mathbb{R}\times\Omega\rightarrow\mathbb{R}$ is defined by
\begin{equation}
J(\lambda,u):=\int_{0}^{2\pi}\left[  \frac{1}{2}|\dot{u}(t)|^{2}-\lambda
^{2}V(u(t))\right]  dt. \label{eq:var-1}%
\end{equation}
Assume that $u_{o}\in\mathbb{R}^{12}$ is the equilibrium point of
\eqref{eq:mol} described in subsection \ref{sec:equilib}. Then clearly,
$u_{o}$ is a critical point of $J$. We are interested in finding
non-stationary $2\pi$-periodic solutions bifurcating from $u_{o}$, i.e.
non-constant solutions to system \eqref{eq:bif1}.

We consider the $G$-orbit of $u_{o}$ in the space $H^{1}(S^{1},\mathscr V)$.
We denote by $\mathcal{S}_{o}$ the slice to $G(u_{o})$ in $H^{1}%
(S^{1},\mathscr V)$. We will also denote by
\[
\mathscr J:\mathbb{R}\times\left(  \mathcal{S}_{o}\cap\Omega\right)
\rightarrow\mathbb{R}%
\]
the restriction of $J$ to the set $\mathcal{S}_{o}\cap\Omega$. Then clearly,
$\mathscr J$ is $G_{u_{o}}$-invariant. Then, since the orbit $G(u_{o})$ is
orthogonal to the slice $\mathcal{S}_{o}$, in a small tubular neighborhood of
the orbit $G(u_{o})$, critical points of $\mathscr J$ are critical points of
$J$ and consequently, they are solutions to system \eqref{eq:bif1}. This
property allows to establish the Slice Criticality Principle (see Theorem
\ref{thm:SCP}), to compute the $G$-equivariant gradient degree of $J$ on this
small tubular neighborhood, which will provide us the full equivariant
topological classification of all non-constant periodic orbits bifurcating
from the equilibrium $u_{o}$.

Consider the operator $L:H^{2}(S^{1};\mathscr V)\rightarrow L^{2}%
(S^{1};\mathscr V)$, given by $Lu=-\ddot{u}+u$, $u\in H^{2}(S^{1},\mathscr
V)$. Then the inverse operator $L^{-1}$ exists and is bounded. Put
$j:H^{2}(S^{1};\mathscr V)\rightarrow H^{1}(S^{1},\mathscr V)$ be the natural
embedding operator. Clearly, $j$ is a compact operator. Then, one can easily
verify that
\begin{equation}
\nabla_{u}J(\lambda,u)=u-j\circ L^{-1}(\lambda^{2}\nabla V(u)+u),
\label{eq:gradJ}%
\end{equation}
where $u\in H^{1}(S^{1},\mathscr V)$. Consequently, the bifurcation problem
\eqref{eq:bif1} can be written as $u-j\circ L^{-1}(\lambda^{2}\nabla
V(u)+u)=0$. Moreover, we have
\begin{equation}
\nabla_{u}^{2}J(\lambda,u_{o})v=v-j\circ L^{-1}(\lambda^{2}\nabla^{2}%
V(u_{o})v+v)~, \label{eq:D2J}%
\end{equation}
where $v\in H^{1}(S^{1},\mathscr V)$.

Consider the operator
\begin{equation}
\mathscr A(\lambda):=\nabla_{u}^{2}J(\lambda,u_{o})|_{\mathcal{S}_{o}%
}:\mathcal{S}_{o}\rightarrow\mathcal{S}_{o}. \label{eq:opA}%
\end{equation}
Notice that
\[
\nabla_{u}^{2}\mathscr J(\lambda,u_{o})=\mathscr A(\lambda),
\]
thus, by implicit function theorem, $G(u_{o})$ is an isolated orbit of
critical points of $J$, whenever $\mathscr A(\lambda)$ is an isomorphism.
Therefore, if a point $(\lambda_{o},u_{o})$ is a bifurcation point for
\eqref{eq:bif1}, then $\mathscr A(\lambda_{o})$ cannot be an isomorphism. In
such a case we put
\[
\Lambda:=\{\lambda>0:\mathscr A(\lambda_{o})\text{ is not an isomorphism}\}~,
\]
and will call the set $\Lambda$ the \textit{critical set} for the trivial
solution $u_{o}$.

\subsection{Bifurcation Theorem}

Consider the $S^{1}$-action on $H^{1}(S^{1},\mathscr V)$, where $S^{1}$ acts
on functions by shifting the argument (see \eqref{eq:ac}). Then, $(H^{1}%
(S^{1},\mathscr V))^{S^{1}}$ is the space of constant functions, which can be
identified with the space $\mathscr V$, i.e.
\[
H^{1}(S^{1},\mathscr V)=\mathscr V\oplus\mathscr W,\quad\mathscr W:=\mathscr
V^{\perp}.
\]
Then, the slice $\mathcal{S}_{o}$ in $H^{1}(S^{1},\mathscr V)$ to the orbit
$G(u_{o})$ at $u_{o}$ is exactly
\[
\mathcal{S}_{o}=S_{o}\oplus\mathscr W.
\]

Any $\lambda_{o}\in\Lambda$ satisfies the condition that $\ \mathscr
A(\lambda_{o})|_{S_{o}}:S_{o}\rightarrow S_{o}$ is an isomorphism, since the
eigenvalues are $\mu_{j}\neq0$ for $j=0,1,2$, which leads to:

\begin{theorem}
\label{th:bif1} Consider the bifurcation system \eqref{eq:bif1} and assume
that $\lambda_{o}\in\Lambda$ is isolated in the critical set $\Lambda$, i.e.
there exists $\lambda_{-}<\lambda_{o}<\lambda_{+}$ such that $[\lambda
_{-},\lambda_{+}]\cap\Lambda=\{\lambda_{o}\}$. Define
\[
\omega_{G}(\lambda_{o}):=\nabla_{G_{u_{o}}}\text{\textrm{-deg}}\Big(\mathscr
A(\lambda_{-}),B_{1}(0)\Big)-\nabla_{G_{u_{o}}}\text{\textrm{-deg}%
}\Big(\mathscr A(\lambda_{+}),B_{1}(0)\Big),
\]
where $B_{1}(0)$ stands for the open unit ball in $\mathscr H$. If
\[
\omega_{G}(\lambda_{o})=n_{1}(H_{1})+n_{2}(H_{2})+\dots+n_{m}(H_{m})
\]
is non-zero, i.e. $n_{j}\not =0$ for some $j=1,2,\dots,m$, then there exists a
bifurcating branch of nontrivial solutions to \eqref{eq:bif1} from the orbit
$\{\lambda_{o}\}\times G(u_{o})$ with symmetries at least $(H_{j})$.
\end{theorem}

Consider the $S^{1}$-isotypic decomposition of $\mathscr W$, i.e.
\[
\mathscr W=\overline{\bigoplus_{l=1}^{\infty}\mathscr W_{l}},\quad\mathscr
W_{l}:=\{\cos(l\cdot)\mathfrak{a}+\sin(l\cdot)\mathfrak{b}:\mathfrak{a}%
,\,\mathfrak{b}\in\mathscr V\}
\]
In a standard way, the space $\mathscr W_{l}$, $l=1,2,\dots$, can be naturally
identified with the space $\mathscr V^{\mathbb{C}}$ on which $S^{1}$ acts by
$l$-folding,
\[
\mathscr W_{l}=\{e^{il\cdot}z:z\in\mathscr V^{\mathbb{C}}\}.
\]

Since the operator $\mathscr A(\lambda)$ is $G_{u_{o}}$-equivariant with
\[
G_{u_{o}}=\tilde{S}_{4}\times O(2),
\]
it is also $S^{1}$-equivariant and thus $\mathscr A(\lambda)(\mathscr
W_{l})\subset\mathscr W_{l}$. Using the $\tilde{S}_{4}$-isotypic decomposition
of $\mathscr V^{\mathbb{C}}$, we have the $G_{u_{o}}$-isotypic decomposition
\[
\mathscr W_{l}=W_{0,l}\oplus W_{1,l}\oplus W_{2,l}\oplus W_{3,l}~,\qquad
W_{j,l}=\mathcal{W}_{j,l}~.
\]
Moreover, we have
\[
\mathscr A(\lambda)|_{W_{j,l}}=\left(  1-\frac{\lambda^{2}\mu_{j}+1}{l^{2}%
+1}\right)  \text{\textrm{Id\,}}~,
\]
which implies that $\lambda_{o}\in\Lambda$ if and only if $\lambda_{o}%
^{2}=l^{2}/\mu_{j}$ for some $l=1,2,3,\dots$ and $j=0,1,2$.

Then the critical set $\Lambda$ for the equilibrium $u_{o}$ of the system
\eqref{eq:mol} is
\[
\Lambda:=\left\{  \frac{l}{\sqrt{\mu_{j}}}:j=0,1,2,\quad l=1,2,3,\dots
\right\}  ,
\]
and we can identify the critical numbers $\lambda\in\Lambda$ as
\[
\lambda_{j,l}=\frac{l}{\sqrt{\mu_{j}}}~.
\]
The critical numbers are not uniquely identified by the indices $(j,l)$ due to
resonances. Indeed, let us list the first critical numbers from $\Lambda$
\[
\lambda_{0,1}<\lambda_{1,1}<\lambda_{2,1}=\lambda_{0,2}<\lambda_{1,2}%
<\lambda_{2,2}=\lambda_{0,4}~.
\]

\begin{definition}
For simplicity, hereafter we denote by $S_{4}$ the isotropy group
$\mathfrak{G}_{u_{0}}=\tilde{S}_{4}$, i.e. with this notation we have that
\[
G_{u_{0}}=\mathfrak{G}_{u_{0}}\times O(2)=S_{4}\times O(2)~.
\]

\end{definition}

From the computation of the gradient degree in (\ref{eq:lin-GdegGrad}) with
$G_{u_{0}}$, we obtain for $\lambda\notin\Lambda$ that
\begin{equation}
\nabla_{G_{u_{0}}}\text{\textrm{-deg}}\Big(\mathscr A(\lambda_{o}%
),B_{1}(0)\Big)=\prod_{\left\{  \left(  j,l\right)  \in\mathbb{N}^{2}%
:\lambda_{j,l}<\lambda_{o}\right\}  }\nabla\text{\textrm{-deg}}_{\mathcal{W}%
_{j,l}}~.\label{Pre}%
\end{equation}

For each critical number $\lambda_{j,l}$ we choose two numbers $\lambda
_{-}<\lambda_{j,l}<\lambda_{+}$ such that $[\lambda_{-},\lambda_{+}%
]\cap\Lambda=\{\lambda_{j,l}\}$. Calculating the difference of the gradient
degree at $\lambda_{+}$ and $\lambda_{-}$ using (\ref{Pre}), we obtain that
the equivariant invariants are given by
\begin{align*}
\omega_{G}(\lambda_{0,1}) &  =\nabla\text{\textrm{-deg}}_{\mathcal{W}_{0,1}%
}-(S_{4}\times O(2))\\
\omega_{G}(\lambda_{1,1}) &  =\nabla\text{\textrm{-deg}}_{\mathcal{W}_{0,1}%
}\ast\Big(\nabla\text{\textrm{-deg}}_{\mathcal{W}_{1,1}}-(S_{4}\times
O(2))\Big)\\
\omega_{G}(\lambda_{2,1}) &  =\nabla\text{\textrm{-deg}}_{\mathcal{W}_{0,1}%
}\ast\nabla\text{\textrm{-deg}}_{\mathcal{W}_{1,1}}\ast\Big(\nabla
\text{\textrm{-deg}}_{\mathcal{W}_{2,1}}\ast\nabla\text{\textrm{-deg}%
}_{\mathcal{W}_{0,2}}-(S_{4}\times O(2))\Big)
\end{align*}

\subsection{Computation of the Gradient Degree}

We consider the product group $G_{1}\times G_{2}$ given two groups $G_{1}$ and
$G_{2}$. The well-known result (see \cite{DKY,Goursat}) provides a description
of the product group $G_{1}\times G_{2}$. Namely, for any subgroup
$\mathscr H$ of the product group $G_{1}\times G_{2}$ there exist subgroups
$H\leq G_{1}$ and $K\leq G_{2}$ a group $L$ and two epimorphisms
$\varphi:H\rightarrow L$ and $\psi:K\rightarrow L$ such that
\begin{equation}
\mathscr H=\{(h,k)\in H\times K:\varphi(h)=\psi(k)\},
\end{equation}
In this case, we will use the notation
\[
\mathscr H=:H\prescript{\varphi}{}\times_{L}^{\psi}K,
\]
and the group $:H\prescript{\varphi}{}\times_{L}^{\psi}K$ will be called an
\textit{amalgamated} subgroup of $G_{1}\times G_{2}$.

Therefore, any closed subgroup $\mathscr H$ of $S_{4}\times O(2)$ is an
amalgamated subgroup $H^{\varphi}\times_{L}^{\psi}K$, where $H\leq S_{4}$ and
$K\leq O(2)$. In order to make amalgamated subgroup notation simpler and
self-contained we will assume that%
\[
L=K/\ker(\psi),
\]
so $\psi:K\rightarrow L$ is evidently the natural projection and there is no
need to indicate it. On the other hand the group $L$ can be naturally
identified with a finite subgroup of $O(2)$ being either $D_{n}$ or
$\mathbb{Z}_{n}$, $n\geq1$. Since we are interested in describing conjugacy
classes of $\mathscr H$, we can identify the epimorphism $\varphi:H\rightarrow
L$ by indicating
\[
Z=\text{Ker\thinspace}(\varphi)\quad\text{ and }\quad R=\varphi^{-1}(\langle
r\rangle)
\]
where $r$ is the rotation generator in $L$ and $\langle r\rangle$ is the
cyclic subgroup generated by $r$. Then, instead of using the notation
$H^{\varphi}\times_{L}^{\psi}K$ we will write
\begin{equation}
\mathscr H=:H{\prescript{Z}{}\times_{L}^{R}}K~, \label{eq:amalg}%
\end{equation}
where $H$, $Z$ and $R$ are subgroups of $S_{4}$ identified by
\begin{align*}
V_{4}  &  =\{(1),(12)(34),(13)(24),(14)(23)\}~,\\
D_{4}  &  =\{(1),(1324),(12)(34),(1423),(34),(14)(23),(12),(13)(24)\}~,\\
Z_{4}  &  =\{(1),(1324),(12)(34),(1423)\}\,,\\
D_{3}  &  =\{(1),(123),(132),(12),(23),(13)\}~,\\
D_{2}  &  =\{(1),(12)(34),(12),(34)\}~,\\
D_{1}  &  =\{(1),(12)\}~.
\end{align*}
In the case when all the epimorphisms $\varphi$ with the kernel $Z$ are
conjugate, there is no need to use the symbol $R$ in \eqref{eq:amalg}, so we
will simply write $\mathscr H=H{\prescript{Z}{}\times_{L}K}$. In addition, in
the case all epimorphisms $\varphi$ from $H$ to $L$ are conjugate, we can also
omit the symbol $Z$, i.e. we will write $\mathscr H=H\times_{L}K$.

The notation explained in this section is useful to obtain the classification
of the all conjugacy classes $(\mathscr H)$ of closed subgroups in
$S_{4}\times O(2)$.

Let us point out that to obtain a complete equivariant classification of the
bifurcating branches of nontrivial solutions, the full topological invariant
$\omega_{G}(\lambda_{j_{o},1})\in U(I\times O(2))$ should be considered. In
particular, although it is not the case here, the invariant $\omega
_{G}(\lambda_{j_{o},1})$ may contain maximal orbit types $(H)$ with infinite
Weyl's group $W(H)$. With the use of GAP programming (see \cite{Pin}) one will
definitely be able to establish the exact value of the invariant $\omega
_{G}(\lambda_{j_{o},1})$, but at this moment requires additional computer
programming. Therefore, in order to simplify the computations, we consider its
truncation to $A(I\times O(2))$, given by
\[
\widetilde{\omega}_{G}(\lambda_{j_{o},1}):=\pi_{0}\Big(\omega_{G}%
(\lambda_{j_{o},1})\Big).
\]
where $\pi_{0}:U(G)\rightarrow A(G)$ is a ring homomorphism. Other more
complex molecular structures may require the full value of the invariant
$\omega_{G}(\lambda_{j_{o}})$ in the Euler ring $U(G)$, we should keep in mind
that it is necessary to use the full $G$-equivariant gradient degree for its analysis.

We can use GAP programming (see \cite{Pin}) to compute the basic degrees
truncated to $A(G)$.
\begin{align*}
\mathrm{Deg}_{\mathcal{W}_{0,l}}=\;  &  -{\color{red}({S_{4}}%
\prescript{}{}\times_{{}}D_{l})}+{({S_{4}}\prescript{}{}\times_{{}}O(2))},\\
\mathrm{Deg}_{\mathcal{W}_{1,l}}=\;  &  -{\color{red}({D_{4}}%
\prescript{{D_2}}{}\times_{\mathbb{Z}_{2}}D_{2l})}-{\color{red}({D_{2}%
}\prescript{{D_1}}{}\times_{\mathbb{Z}_{2}}D_{2l})}-{\color{red}({D_{4}%
}\prescript{{\bz_1}}{}\times_{D_{4}}D_{4l})}-{\color{red}({D_{3}%
}\prescript{}{}\times_{{}}D_{l})}\\
&  -{\color{red}({D_{3}}\prescript{{\bz_1}}{}\times_{D_{3}}D_{3l})}+2{({D_{1}%
}\prescript{}{}\times_{{}}D_{l})}-{({\mathbb{Z}_{1}}\prescript{}{}\times_{{}%
}D_{l})}+{({\mathbb{Z}_{2}}\prescript{{\bz_1}}{}\times_{\mathbb{Z}_{2}}%
D_{2l})}\\
&  +{({D_{2}}\prescript{{\bz_1}}{{\bz_2}}\times_{D_{2}}D_{2l})}+{({V_{4}%
}\prescript{{\bz_1}}{}\times_{D_{2}}D_{2l})}+{({D_{2}}%
\prescript{{D_1}}{}\times_{D_{1}}D_{l})}-{({\mathbb{Z}_{2}}%
\prescript{{\bz_1}}{}\times_{D_{1}}D_{l})}\\
&  +{({S_{4}}\prescript{}{}\times_{{}}O(2))},\\
\mathrm{Deg}_{\mathcal{W}_{2,l}}=\;  &  -{\color{red}({S_{4}}%
\prescript{{V_4}}{}\times_{D_{3}}D_{3l})}-{({D_{4}}\prescript{}{}\times_{{}%
}D_{l})}+{({V_{4}}\prescript{}{}\times_{{}}D_{l})}-{({D_{4}}%
\prescript{{V_4}}{}\times_{\mathbb{Z}_{2}}D_{2l})}\\
&  +2{({D_{4}}\prescript{{V_4}}{}\times_{D_{1}}D_{l})}+{({S_{4}}%
\prescript{}{}\times_{{}}O(2))},\\
\mathrm{Deg}_{\mathcal{W}_{3,l}}=\;  &  -{\color{red}({D_{4}}%
\prescript{{\bz_4}}{}\times_{\mathbb{Z}_{2}}D_{2l})}-{\color{red}({D_{4}%
}\prescript{{\bz_1}}{}\times_{D_{4}}D_{4l})}-{\color{red}({D_{2}%
}\prescript{{D_1}}{}\times_{\mathbb{Z}_{2}}D_{2l})}+2{({D_{1}}%
\prescript{{\bz_1}}{}\times_{\mathbb{Z}_{2}}D_{2l})}\\
&  +{({\mathbb{Z}_{2}}\prescript{{\bz_1}}{}\times_{\mathbb{Z}_{2}}D_{2l}%
)}-{({\mathbb{Z}_{1}}\prescript{}{}\times_{{}}D_{l})}-{({D_{3}}%
\prescript{{\bz_3}}{}\times_{\mathbb{Z}_{2}}D_{2l})}-{({D_{3}}%
\prescript{{\bz_1}}{}\times_{D_{3}}D_{3l})}\\
&  +{({D_{2}}\prescript{{\bz_1}}{{D_1}}\times_{D_{2}}D_{2l})}+{({D_{2}%
}\prescript{{\bz_1}}{{\bz_2}}\times_{D_{2}}D_{2l})}+{({V_{4}}%
\prescript{{\bz_1}}{}\times_{D_{2}}D_{2l})}-{({\mathbb{Z}_{2}}%
\prescript{{\bz_1}}{}\times_{D_{1}}D_{l})}\\
&  +{({S_{4}}\prescript{}{}\times_{{}}O(2))},\\
\mathrm{Deg}_{\mathcal{W}_{4,l}}=\;  &  -{\color{red}({S_{4}}%
\prescript{{A_4}}{}\times_{\mathbb{Z}_{2}}D_{2l})}+{({S_{4}}%
\prescript{}{}\times_{{}}O(2))},
\end{align*}

Next, we use GAP programming (see \cite{Pin}) and the product $\ast$ of the
Euler ring $U(\Gamma)$ to compute the full equivariant invariants to
$A(I\times O(2))$, where the maximal isotropy classes are colored red:%
\begin{align*}
\widetilde{\omega}_{G}(\lambda_{0,1})=\;  &  -{\color{red}({S_{4}%
}\prescript{}{}\times_{{}}D_{1})},\\
\widetilde{\omega}_{G}(\lambda_{1,1})=\;  &  -{\color{red}({D_{4}%
}\prescript{{D_2}}{}\times_{\mathbb{Z}_{2}}D_{2})}-{\color{red}({D_{2}%
}\prescript{{D_1}}{}\times_{\mathbb{Z}_{2}}D_{2})}-{\color{red}({D_{4}%
}\prescript{{\bz_1}}{}\times_{D_{4}}D_{4})}+{\color{red}({D_{3}}%
\prescript{}{}\times_{{}}D_{1})}\\
&  -{\color{red}({D_{3}}\prescript{{\bz_1}}{}\times_{D_{3}}D_{3})}+{({D_{2}%
}\prescript{}{}\times_{{}}D_{1})}-{({D_{1}}\prescript{}{}\times_{{}}D_{1}%
)}+{({\mathbb{Z}_{2}}\prescript{{\bz_1}}{}\times_{\mathbb{Z}_{2}}D_{2})}\\
&  +{({D_{2}}\prescript{{\bz_1}}{{\bz_2}}\times_{D_{2}}D_{2})}+{({V_{4}%
}\prescript{{\bz_1}}{}\times_{D_{2}}D_{2})}+{({D_{4}}\prescript{{D_2}}{}\times
_{D_{1}}D_{1})}+{({D_{1}}\prescript{{\bz_1}}{}\times_{D_{1}}D_{1})}\\
&  -{({\mathbb{Z}_{2}}\prescript{{\bz_1}}{}\times_{D_{1}}D_{1})},\\
\widetilde{\omega}_{G}(\lambda_{2,1})=  &  -{\color{red}({S_{4}}%
\prescript{{V_4}}{}\times_{D_{3}}D_{3})}-{\color{red}({S_{4}}%
\prescript{}{}\times_{{}}D_{2})}+{({D_{4}}\prescript{{V_4}}{}\times
_{\mathbb{Z}_{2}}D_{2})}-{({\mathbb{Z}_{4}}\prescript{{\bz_2}}{}\times
_{\mathbb{Z}_{2}}D_{2})}\\
&  +2{({D_{2}}\prescript{{D_1}}{}\times_{\mathbb{Z}_{2}}D_{2})}-{({D_{1}%
}\prescript{{\bz_1}}{}\times_{\mathbb{Z}_{2}}D_{2})}-2{({\mathbb{Z}_{2}%
}\prescript{{\bz_1}}{}\times_{\mathbb{Z}_{2}}D_{2})}+2{({S_{4}}%
\prescript{}{}\times_{{}}D_{1})}\\
&  -{({D_{4}}\prescript{}{}\times_{{}}D_{1})}-2{({D_{3}}\prescript{}{}\times
_{{}}D_{1})}-{({D_{2}}\prescript{}{}\times_{{}}D_{1})}+{({D_{1}}%
\prescript{}{}\times_{{}}D_{1})}\\
&  +{({\mathbb{Z}_{1}}\prescript{}{}\times_{{}}D_{1})}+2{({D_{4}%
}\prescript{{D_2}}{}\times_{\mathbb{Z}_{2}}D_{2})}+2{({D_{3}}%
\prescript{{\bz_1}}{}\times_{D_{3}}D_{3})}-{({D_{4}}%
\prescript{{\bz_2}}{{\bz_4}}\times_{D_{2}}D_{2})}\\
&  -{({D_{2}}\prescript{{\bz_1}}{{D_1}}\times_{D_{2}}D_{2})}-{({D_{2}%
}\prescript{{\bz_1}}{{\bz_2}}\times_{D_{2}}D_{2})}-{({V_{4}}%
\prescript{{\bz_1}}{}\times_{D_{2}}D_{2})}-{({D_{4}}\prescript{{D_2}}{}\times
_{D_{1}}D_{1})}\\
&  +{({D_{4}}\prescript{{V_4}}{}\times_{D_{1}}D_{1})}-{({D_{2}}%
\prescript{{D_1}}{}\times_{D_{1}}D_{1})}+3{({\mathbb{Z}_{2}}%
\prescript{{\bz_1}}{}\times_{D_{1}}D_{1})}.
\end{align*}

\section{Description of the Symmetries}

The invariants $\omega_{G}(\lambda_{j,1})$ give the bifurcation of periodic
solutions for each of five maximal groups. However, we only know that a group
is maximal if it is maximal in a certain isotypical component of a Fourier
mode. Since the bifurcation from $\lambda_{2,1}=\lambda_{0,2}$ with maximal
group ${S_{4}\prescript{S_4}{}\times_{\mathbb{Z}_{1}}^{{}}D_{2}}$ is not
independent of minimal period bifurcation from $\lambda_{0,1}$ with maximal
group ${S_{4}\prescript{S_4}{}\times_{\mathbb{Z}_{1}}^{{}}D_{1}}$, we cannot
conclude that these two bifurcation are different from each other.

We can conclude that the other 7 maximal groups in the invariants $\omega
_{G}(\lambda_{j,1})$ for $j=0,1,2$ give different global families of periodic
solutions with period $T=2\pi\lambda_{j,1}l_{o}$ (for some $l_{o}\in
\mathbb{N}$), where $(\lambda_{j,1}l_{o})^{-1}$ is the limit frequency. Next
we describe the symmetries of the solutions for these maximal isotropy groups.
Notice that we have identified the elements of $\tilde{S}_{4}$ with $S_{4}$,
i.e. an element $\sigma\in S_{4}$ in a maximal group acts as%
\[
\sigma u_{j}=A_{\sigma}u_{\sigma(j)}~.
\]

\subsection{Families with Frequency $\sqrt{\mu_{0}}$}

The tetrahedron configuration has one global family of periodic solutions
starting with frequency $\lambda_{0,1}^{-1}=\sqrt{\mu_{0}}$. This family has
symmetries
\[
S_{4}\prescript{S_4}{}\times_{\mathbb{Z}_{1}}^{{}}D_{1}.
\]
This group is generated by $S_{4}$ and $\kappa\in D_{1}$. The symmetry $S_{4}$
implies that the configurations is a regular tetrahedron at any time.
Moreover, the group $D_{1}$ implies that%
\[
u(t)=\kappa u(t)=u(-t),
\]
i.e. the periodic solution is a brake orbit, which means that the velocity
$\dot{u}$ of all the molecules are zero at the times $t=0,\pi$,
\[
\dot{u}(0)=\dot{u}(\pi)=0.
\]
Therefore, these solution consist of a regular tetrahedron that expands and
contracts in periodic motion, in an orbit which is similar to a line.

\subsection{Families with Frequency $\sqrt{\mu_{1}}$}

The tetrahedron configuration has five different families of periodic
solutions starting with frequency $\lambda_{1,1}^{-1}=\sqrt{\mu_{1}}$, each
family with a different group of symmetries.

The group
\[
{D_{4}\prescript{D_2}{}\times_{\mathbb{Z}_{2}}^{{}}D_{2}}%
\]
is generated by the elements $\kappa\in O(2)$, $(12),(34)\in S_{4}$ and
$((13)(24),e^{i\pi})\in S_{4}\times O(2)$. The element $\kappa$ implies that
the periodic solution is a brake orbit as in the description before. We have
that $A_{(12)}$ is the inversion over the plane that has the points
$\gamma_{3}$, $\gamma_{4}$ and the middle point of $\gamma_{1}$ and
$\gamma_{2}$, and such that it interchanges $\gamma_{1}$ with $\gamma_{2}$.
Then, the symmetry $(12)$ implies that $u_{1}$ is the inversion of $u_{2}$,
and similarly $(34)$ implies that $u_{3}$ is the inversion of $u_{4}$. The
element $A_{(13)(24)}$ is a rotation by $\pi$ that interchanges $\gamma_{1}$
with $\gamma_{3}$ and $\gamma_{2}$ with $\gamma_{4}$. Therefore $u_{1}$ is the
$\pi$-rotation and $\pi$-phase shift of $u_{3}(t)$. In this symmetries all the
orbits are determined by the positions of only one of the particles $u_{1}$.

The group
\[
{D_{2}\prescript{D_1}{}\times_{\mathbb{Z}_{2}}^{{}}D_{2}}%
\]
is generated by the elements $\kappa\in O(2)$, $(12)\in S_{4}$ and
$((34),e^{i\pi})\in S_{4}\times O(2)$. The symmetry $\kappa$ implies that the
solutions is a brake orbit, $(12)$ that $u_{1}$ is the inversion of $u_{2}$,
and $((34),e^{i\pi})$ that $u_{3}$ is the $\pi$-rotation and $\pi$-phase shift
of $u_{4}(t)$. In this case, the orbit of $u_{1}$ determines $u_{2}$ and
$u_{4}$ determines $u_{3}$, but there is no relation among these two pairs of particles.

The group
\[
{D_{4}\prescript{\mathbb{Z}_1}{}\times_{D_{4}}^{{}}D_{4}}%
\]
is generated by $\left(  (12),\kappa\right)  $ and $\left(  (1324),e^{i\pi
/2}\right)  $ in $S_{4}\times O(2)$. The element $\left(  (12),\kappa\right)
$ implies that $u_{1}(t)$ is the inversion of $u_{2}(-t)$. In this case, the
orbit is not brake, which means it is similar to a circle. The matrix
$A_{(1324)}$ is a rotor reflections by $\pi/2$. Then the symmetry $\left(
(1324),e^{i\pi/2}\right)  $ implies that the particles are related by applying
a $\pi/2$-rotoreflection and at the same, a temporal phase shift by $\pi/2$.

The group
\[
{{D_{3}}\prescript{}{}\times_{{}}D_{1}}%
\]
is generated by $\kappa$, which implies that the solution is a brake orbits,
and the group $D_{3}$, which implies that the positions $u_{1}$, $u_{2}$ and
$u_{3}$ always form a triangle. In this case, the position $u_{4}$ follows a
movement that counterbalance the triangle formed by these elements.

The group%
\[
-{{D_{3}}\prescript{{\bz_1}}{}\times_{D_{3}}D_{3}}%
\]
is generated by the elements $((123),e^{i2\pi/3})$ and $((12),\kappa)$. The
element\break$((123),e^{i2\pi/3})$ implies that $u_{1}(t)=u_{2}(t+2\pi
/3)=u_{3}(t+4\pi/3)$ and therefore, the movement in these three elements is a
(discrete)\ rotating wave. In addition, the element $\left(  (12),\kappa
\right)  $ implies that this rotating wave is invariant by an inversion in
time, $u_{1}(t)=A_{(1,2)}u_{2}(-t)=A_{(1,2)}u_{1}(-t-2\pi/3)$.

\subsection{Families with Fequency $\sqrt{\mu_{2}}$}

The tetrahedron configuration has one family of periodic solutions starting
with frequency $\lambda_{1,1}^{-1}=\sqrt{\mu_{2}}$ with symmetries
${S_{4}\prescript{V_4}{}\times_{D_{3}}^{{}}D_{3}}$. This group is generated by
$V_{4}$ and $\left(  (123),2\pi/3\right)  ,((12),\kappa)\in S_{4}\times O(2)$.
The symmetries $V_{4}$ place coordinates $u(t)$ in form of a non-regular
tetrahedron figure with two axis of symmetry at any time. The element $\left(
(123),2\pi/3\right)  $ means that $u_{1}$, $u_{2}$ and $u_{3}$ are related by
a rotation of $2\pi/3$ and a phase shift of $2\pi/3$.

\section{Appendix: \textbf{Equivariant Gradient Degree}}

\subsection{Group Actions}

In what follows $G$ always stands for a compact Lie group and all subgroups of
$G$ are assumed to be closed \cite{Bre, Kawa}. For a subgroup $H\subset G$,
denote by $N\left(  H\right)  $ the \textit{normalizer} of $H$ in $G$, and by
$W\left(  H\right)  =N\left(  H\right)  /H$ the \textit{Weyl group} of $H$ in
$G$. In the case when we are dealing with different Lie groups, we also write
$N_{G}\left(  H\right)  $ (resp. $W_{G}\left(  H\right)  $) instead of
$N\left(  H\right)  $ (resp. $W\left(  H\right)  $). We denote by $\left(
H\right)  $ the conjugacy class of $H$ in $G$ and use the notations:
\begin{align*}
\Phi\left(  G\right)   &  :=\left\{  \left(  H\right)  :H\;\;\text{is a
subgroup of }\;G\right\}  ,\\
\Phi_{n}\left(  G\right)   &  :=\left\{  \left(  H\right)  \in\Phi\left(
G\right)  :\text{\textrm{dim\thinspace}}W\left(  H\right)  =n\right\}  .
\end{align*}
The set $\Phi\left(  G\right)  $ has a natural partial order defined by
\begin{equation}
\left(  H\right)  \leq\left(  K\right)  \;\;\Longleftrightarrow\;\;\exists
_{g\in G}\;\;gHg^{-1}\subset K. \label{eq:partial}%
\end{equation}

For a $G$-space $X$ and $x\in X$, the subgroup $G_{x}:=\left\{  g\in
G:\;gx=x\right\}  $ is called the \textit{isotropy} of $x$; $G\left(
x\right)  :=\left\{  gx:\;g\in G\right\}  $ the \textit{orbit} of $x$, and the
conjugacy class $\left(  G_{x}\right)  $ is called the \textit{orbit type} of
$x$. Also, for a subgroup $H\subset G$, we use
\[
X^{H}:=\left\{  x\in X:\;G_{x}\supset H\right\}
\]
for the fixed point space of $H$. The orbit space for a $G$-space $X$ will be
denoted by $X/G$.

As any compact Lie group admits only countably many non-equivalent real (resp.
complex) irreducible representations. Given a compact Lie group $G$, we will
assume that we know a complete list of all its real (resp. complex)
irreducible representations, denoted $\mathcal{V}_{i}$, $i=0,$ $1,$ $\ldots$
(resp. $\mathcal{W}_{j}$, $j=0,$ $1,$ $\ldots$). We refer to \cite{AED} for
examples of such lists and the related notation.

Let $V$ (resp. $W$) be a finite-dimensional real (resp. complex) $\Gamma
$-representation (without loss of generality, $V$ (resp. $W$) can be assumed
to be orthogonal (resp. unitary)). Then, $V$ (resp. $W$) decomposes into the
direct sum of $G$-invariant subspaces
\begin{equation}
V=V_{0}\oplus V_{1}\oplus\dots\oplus V_{r}\text{,} \label{eq:Giso}%
\end{equation}
(resp. $W=W_{0}\oplus W_{1}\oplus\dots\oplus W_{s}$), called the $G$%
\textit{-}\emph{isotypical decomposition of }$V$ (resp. $W$), where each
isotypical component $V_{i}$ (resp. $W_{j}$) is \emph{modeled} on the
irreducible $G$-representation $\mathcal{V}_{i}$, $i=0,$ $1,$ $\dots,$ $r$,
(resp. $\mathcal{W}_{j}$, $j=0,$ $1,$ $\dots,$ $s$), i.e. $V_{i}$ (resp.
$W_{j}$) contains all the irreducible subrepresentations of $V$ (resp. $W$)
which are equivalent to $\mathcal{V}_{i}$ (resp. $\mathcal{W}_{j}$).

\subsection{Euler Ring}

Let
\[
U\left(  G\right)  :={\mathbb{Z}}\left[  \Phi\left(  G\right)  \right]  .
\]
denotes the free $Z$-module generated by $\Phi(G)$ .

\begin{definition}
\textrm{\textrm{\label{def:EulerRing}$($cf. $\mathrm{\cite{tD}})$ Define a
ring multiplication on generators $\left(  H\right)  $, $\left(  K\right)
\in\Phi\left(  G\right)  $ as follows:
\begin{equation}
\left(  H\right)  \ast\left(  K\right)  =\sum_{\left(  L\right)  \in
\Phi\left(  G\right)  }n_{L}\left(  L\right)  ,\label{eq:Euler-mult}%
\end{equation}
where
\begin{equation}
n_{L}:=\chi_{c}\left(  \left(  G/H\times G/K\right)  _{L}/N\left(  L\right)
\right) \label{eq:Euler-coeff}%
\end{equation}
with $\chi_{c}$ the Euler characteristic taken in Alexander-Spanier cohomology
with compact support (cf. \cite{Spa}). The $\mathbb{Z}$-module $U\left(
G\right)  $ equipped with the multiplication \eqref{eq:Euler-mult},
\eqref{eq:Euler-coeff} is a ring called the \textit{Euler ring} of the group
$G$ (cf. \cite{BtD}) } }
\end{definition}

The ${\mathbb{Z}}$-module $A\left(  G\right)  :={\mathbb{Z}}\left[  \Phi
_{0}\left(  G\right)  \right]  $ equipped with a similar multiplication as in
$U\left(  G\right)  $ but restricted to generators from $\Phi_{0}\left(
G\right)  $, is called a \emph{Burnside ring}, i.e.
\[
\left(  H\right)  \cdot\left(  K\right)  =\sum_{\left(  L\right)  }%
n_{L}\left(  L\right)  ,\qquad\left(  H\right)  ,\text{ }\left(  K\right)
,\text{ }\left(  L\right)  \in\Phi_{0}\left(  G\right)  ,
\]
where $n_{L}:=\left(  \left(  G/H\times G/K\right)  _{L}/N\left(  L\right)
\right)  =\left\vert \left(  G/H\times G/K\right)  _{L}/N\left(  L\right)
\right\vert $ (here $\chi$ stands for the usual Euler characteristic). In this
case, we have
\begin{equation}
n_{L}=\frac{n\left(  L,K\right)  \left\vert W\left(  K\right)  \right\vert
n\left(  L,\text{ }H\right)  \left\vert W\left(  H\right)  \right\vert
-\sum_{\left(  I\right)  >\left(  L\right)  }n\left(  L,I\right)
n_{I}\left\vert W\left(  I\right)  \right\vert }{\left\vert W\left(  L\right)
\right\vert }, \label{eq:rec-coef}%
\end{equation}
where
\[
n(L,K)=\left\vert \frac{N(L,K)}{N(K)}\right\vert ,\quad N(L,K):=\{g\in
G:gLg^{-1}\subset K\},
\]
and $\left(  H\right)  ,$ $\left(  K\right)  ,$ $\left(  L\right)  ,$ $\left(
I\right)  $ are taken from $\Phi_{0}\left(  G\right)  $.

Notice that $A\left(  G\right)  $ is a ${\mathbb{Z}}$-submodule of $U\left(
G\right)  $, but not a subring. Define $\pi_{0}:U\left(  G\right)  \rightarrow
A\left(  G\right)  $ on generators $\left(  H\right)  \in\Phi\left(  G\right)
$ by
\begin{equation}
\pi_{0}\left(  \left(  H\right)  \right)  =%
\begin{cases}
\left(  H\right)  & \text{ if }\;\left(  H\right)  \in\Phi_{0}\left(
G\right)  ,\\
0 & \text{ otherwise.}%
\end{cases}
\label{eq:pi_0-homomorphism}%
\end{equation}
Then we have:

\begin{lemma}
\label{lem:pi_0-homomorphism} (cf. \cite{BKR}) The map $\pi_{0}$ defined by
$(\mathrm{\ref{eq:pi_0-homomorphism}})$ is a ring homomorphism, i.e.
\[
\pi_{0}\left(  \left(  H\right)  \ast\left(  K\right)  \right)  =\pi
_{0}\left(  \left(  H\right)  \right)  \cdot\pi_{0}\left(  \left(  K\right)
\right)  ,\qquad\left(  H\right)  ,\text{ }\left(  K\right)  \in\Phi\left(
G\right)  .
\]

\end{lemma}

Let us point out that, although the computations of the Euler ring structure
$U(G)$ are quite challenging in general, in the case $G=\Gamma\times O(2)$
(here $\Gamma$ is a finite group) it can be effectively computed using the
Burnside ring multiplication structure in $A\left(  G\right)  $ via Lemma
\ref{lem:pi_0-homomorphism} allows us to use to partially describe the Euler
ring multiplication structure in $U\left(  G\right)  $.

\subsection{Equivariant Gradient Degree}

Let $G$ be a compact Lie group and $V$ be a $G$-representation. Let
$\varphi:V\rightarrow\mathbb{R}$ be a continuously differentiable
$G$-invariant functional. Define $\mathcal{M}_{\nabla}^{G}$ as the set of
pairs $(\nabla\varphi,\Omega)$ such that $\nabla\varphi$ is an $G$-equivariant
field $\nabla\varphi:V\rightarrow V$ such that
\[
\nabla\varphi(v)\neq0\text{ for }v\in\partial\Omega.
\]

\begin{theorem}
\label{thm:Ggrad-properties} There exists a unique map $\nabla_{G}%
\text{\textrm{-deg\thinspace}}:\mathcal{M}_{\nabla}^{G}\rightarrow U(G)$,
which assigns to every $(\nabla\varphi,\Omega)\in\mathcal{M}_{\nabla}^{G}$ an
element $\nabla_{G}\text{\textrm{-deg\thinspace}}(\nabla\varphi,\Omega)\in
U(G)$, called the $G$\textit{-gradient degree} of $\nabla\varphi$ on $\Omega
$,
\begin{equation}
\nabla_{G}\text{\textrm{-deg\thinspace}}(\nabla\varphi,\Omega)=\sum
_{(H_{i})\in\Phi(\Gamma)}n_{H_{i}}(H_{i})=n_{H_{1}}(H_{1})+\dots+n_{H_{m}%
}(H_{m}), \label{eq:grad-deg}%
\end{equation}
satisfying the following properties:

\begin{description}
\item \textbf{(Existence)} If $\nabla_{G}\text{\textrm{-deg\thinspace}}%
(\nabla\varphi,\Omega)\not =0$, i.e. there is in \eqref{eq:grad-deg} a
non-zero coefficient $n_{H_{i}}$, then exists $u_{0}\in\Omega$ such that
$\nabla\varphi(u_{0})=0$ and $(G_{u_{0}})\geq(H_{i})$.

\item \textbf{(Additivity)} Let $\Omega_{1}$ and $\Omega_{2}$ be two disjoint
open $G$-invariant subsets of $\Omega$ such that $(\nabla\varphi)^{-1}%
(0)\cap\Omega\subset\Omega_{1}\cup\Omega_{2}.$ Then, $\nabla_{G}%
$\textrm{-deg\thinspace}$(\nabla\varphi,\Omega)=\nabla_{G}$%
\textrm{-deg\thinspace}$(\nabla\varphi,\Omega_{1})+\nabla_{G}$%
\textrm{-deg\thinspace}$(\nabla\varphi,\Omega_{2}).$

\item \textbf{(Homotopy)} If $\nabla_{x}\Psi:[0,1]\times V\rightarrow V$ is a
$G$-gradient $\Omega$-admissible homotopy, then
\[
\nabla_{G}\text{\textrm{-deg\thinspace}}(\nabla_{v}\Psi(t,v),\Omega
)=\text{\textit{constant}}.
\]

\item \textbf{(Normalization)} Let $\varphi\in C_{G}^{2}(V,\mathbb{R})$ be a
special $\Omega$-Morse function such that $(\nabla\varphi)^{-1}(0)\cap
\Omega=G(u_{0})$ and $G_{u_{0}}=H$. Then,
\[
\nabla_{G}\text{\textrm{-deg\thinspace}}(\nabla\varphi,\Omega
)=(-1)^{\mathrm{m}^{-}(\nabla^{2}\varphi(u_{0}))}\cdot(H),
\]
where \textquotedblleft$\mathrm{m}^{-}(\cdot)$\textquotedblright\ stands for
the total dimension of eigenspaces for negative eigenvalues of a (symmetric) matrix.

\item \textbf{(Multiplicativity)} For all $(\nabla\varphi_{1},\Omega_{1})$,
$(\nabla\varphi_{2},\Omega_{2})\in\mathcal{M}_{\nabla}^{G}$,
\[
\nabla_{G}\text{\textrm{-deg\thinspace}}(\nabla\varphi_{1}\times\nabla
\varphi_{2},\Omega_{1}\times\Omega_{2})=\nabla_{G}\text{\textrm{-deg\thinspace
}}(\nabla\varphi_{1},\Omega_{1})\ast\nabla_{G}\text{\textrm{-deg\thinspace}%
}(\nabla\varphi_{2},\Omega_{2})
\]
where the multiplication `$\ast$' is taken in the Euler ring $U(G)$.

\item \textbf{(Suspension)} \textit{If }$W$\textit{ is an orthogonal }%
$G$\textit{-representation and $\mathcal{B}$ an open bounded invariant
neighborhood of }$0\in W$\textit{, then}
\[
\nabla_{G}\text{\textrm{-deg\thinspace}}(\nabla\varphi\times\mbox{\rm Id}_{W}%
,\Omega\times\mathcal{B})=\nabla_{G}\text{\textrm{-deg\thinspace}}%
(\nabla\varphi,\Omega).
\]

\item \textbf{(Hopf Property)} Assume $B(V)$ is the unit ball of an orthogonal
$G$-representa\-tion $V$ and for $(\nabla\varphi_{1},B(V)),(\nabla\varphi
_{2},B(V))\in\mathcal{M}_{\nabla}^{G}$, one has
\[
\nabla_{G}\text{\textrm{-deg\thinspace}}(\nabla\varphi_{1},B(V))=\nabla
_{G}\text{\textrm{-deg\thinspace}}(\nabla\varphi_{2},B(V)).
\]
Then, $\nabla\varphi_{1}$ and $\nabla\varphi_{2}$ are $G$-gradient
$B(V)$-admissible homotopic.
\end{description}
\end{theorem}

\subsection{Equivariant Gradient Degree in Hilbert Spaces}

Let $\mathscr H$ be a Hilbert $G$-representation and $\Omega\subset\mathscr H$
and open bounded $G$-invariant set. A $C^{1}$-differentiable $G$-invariant
functional and $f:\mathscr{H}\rightarrow\mathbb{R}$ given by $f(x)=\frac{1}%
{2}\Vert x\Vert^{2}-\varphi(x)$, $x\in\mathscr{H}$, is called \textit{$\Omega
$-admissible} if $\nabla\varphi:\mathscr{H}\rightarrow\mathscr{H}$ is a
completely continuous map and
\[
\forall_{x\in\partial\Omega}\qquad\nabla f(x)=x-\nabla\varphi(x)\not =0~.
\]
By a \textit{$G$-equivariant approximation scheme} $\{P_{n}\}_{n=1}^{\infty}$
in $\mathscr{H}$, we mean a sequence of $G$-equivariant orthogonal projections
$P_{n}:\mathscr{H}\rightarrow\mathscr{H}$, $n=1$, $2$, \dots, such that:

\begin{itemize}
\item[(a)] the subspaces $\mathscr{H}^{n}:=P_{n}(\mathscr{H})$, $n=1,2,$
\dots, are finite-dimensional;

\item[(b)] $\mathscr{H}^{n}\subset\mathscr{H}^{n+1}$, $n=0,1,2,$ \dots;

\item[(c)] $\displaystyle \lim_{n\to\infty} P_{n}x=x$ for all $x\in\mathscr{H}
$.
\end{itemize}

\vskip.3cm

Then for an $\Omega$-admissible $G$-map $f:\mathscr H\rightarrow\mathbb{R}$,
one can define a sequence $f_{n}:\mathscr{H}_{n}\rightarrow\mathbb{R}$ by
$f_{n}(x):=\frac{1}{2}\Vert x\Vert^{2}-\varphi(x)$, $x\in\mathscr{H}_{n}$. By
a standard argument, for sufficiently large $n\in\mathbb{N}$, the maps $\nabla
f_{n}(x):=x-P_{n}\nabla\varphi(x)$, $x\in\mathscr{H}$, are $\Omega_{n}%
$-admissible, where $\Omega_{n}:=\Omega\cap\mathscr{H}_{n}$. Moreover, the
gradient equivariant degrees $\nabla_{G}\text{\textrm{-deg\thinspace}}(\nabla
f_{n},\Omega_{n})$ are well defined and are the same, i.e. for $n$
sufficiently large
\[
\nabla_{G}\text{\textrm{-deg\thinspace}}(\nabla f_{n},\Omega_{n})=\nabla
_{G}\text{\textrm{-deg\thinspace}}(\nabla f_{n+1},\Omega_{n+1}),
\]
which implies by Suspension Property of the $G$-equivariant gradient degree
that we can put
\begin{equation}
\nabla_{G}\text{\textrm{-deg\thinspace}}(\nabla f,\Omega):=\nabla
_{G}\text{\textrm{-deg\thinspace}}(\nabla f_{n},\Omega_{n}), \label{eq:grad-H}%
\end{equation}
$\;\;$ where $\;n\;$ is sufficiently large. One can verify that this
construction doesn't depend on the choice of a $G$-approximation scheme in the
space $\mathscr{H}$, for instance see \cite{DKY}. We should mention that the
ideas behind the usage of the approximation methods to define topological
degree can be rooted to \cite{BP}.

\subsection{Degree on the Slice}

Suppose that the orbit $G(u_{o})$ of $u_{o}\in\mathscr H$ is contained in a
finite-dimensional $G$-invariant subspace, so the $G$-action on that subspace
is smooth and $G(u_{o})$ is a smooth submanifold of $\mathscr H$. Denote by
$S_{o}\subset\mathscr H$ the slice to the orbit $G(u_{o})$ at $u_{o}$. Denote
by $V_{o}:=T_{u_{o}}G(u_{o})$ the tangent space to $G(u_{o})$ at $u_{o}$. Then
clearly, $S_{o}=V_{o}^{\perp}$ and $S_{o}$ is a smooth Hilbert $G_{u_{o}}$-representation.

Then we have (cf. \cite{xx}).

\begin{theorem}
{\smc(Slice Principle)} \label{thm:SCP} Let $\mathscr{E}$ be an orthogonal
$G$-representation, $\varphi:\mathscr{H}\rightarrow\mathbb{R}$ be a
continuously differentiable $G$-invariant functional, $u_{o}\in\mathscr H$ and
$G(u_{o})$ be an isolated critical orbit of $\varphi$ . Let $S_{o}$ be the
slice to the orbit $G(u_{o})$ and $\mathcal{U}$ an isolated tubular
neighborhood of $G(u_{o})$. Put $\varphi_{o}:S_{o}\rightarrow\mathbb{R}$ by
$\varphi_{o}(v):=\varphi(u_{o}+v)$, $v\in S_{o}$. Then
\begin{equation}
\nabla_{G}\text{\textrm{-deg}}(\nabla\varphi,\mathcal{U})=\Theta
(\nabla_{G_{u_{o}}}\text{\textrm{-deg}}(\nabla\varphi_{o},\mathcal{U}\cap
S_{o})), \label{eq:SDP}%
\end{equation}
where $\Theta:U(G_{u_{o}})\rightarrow U(G)$ is defined on generators
$\Theta(H)=(H)$, $(H)\in\Phi(G_{u_{o}})$.
\end{theorem}

We show how to compute $\nabla_{G}\text{-deg\thinspace}(\mathscr A,B(V))$,
where\textrm{ }$\mathscr A:V\rightarrow V$ is a symmetric $G$-equivariant
linear isomorphism and $V$ is an orthogonal $G$-representation, i.e.
$\mathscr
A=\nabla\varphi$ for $\varphi(v)=\frac{1}{2}(\mathscr Av\bullet v)$, $v\in V$,
where \textquotedblleft$\bullet$\textquotedblright\ stands for the inner
product. Consider the $G$-isotypical decomposition \eqref{eq:Giso} of $V$ and
put
\[
\mathscr A_{i}:=\mathscr A|_{V_{i}}:V_{i}\rightarrow V_{i},\quad
i=0,1,\dots,r.
\]
Then, by the multiplicativity property,
\begin{equation}
\nabla_{G}\mbox{-deg}(\mathscr A,B(V))=\prod_{i}^{r}\nabla_{G}%
\mbox{-deg}(\mathscr A_{i},B(V_{i})) \label{eq:deg-Lin-decoGrad}%
\end{equation}
Take $\xi\in\sigma_{-}(\mathscr A)$, where $\sigma_{-}(\mathscr A)$ stands for
the negative spectrum of $\mathscr A$, and consider the corresponding
eigenspace $E(\xi):=\ker(\mathscr A-\xi\mbox{Id})$. Define the numbers
$m_{i}(\xi)$ by
\begin{equation}
m_{i}(\xi):=\dim\left(  E(\xi)\cap V_{i}\right)  /\dim\mathcal{V}_{i},
\label{eq:m_j(mu)-gra}%
\end{equation}
and the so-called basic gradient degrees by
\begin{equation}
\text{deg}_{\mathcal{V}_{i}}:=\nabla_{G}\mbox{-deg}(-\mbox{Id\,},B(\mathcal{V}%
_{i})). \label{eq:basicGrad-deg0}%
\end{equation}
Then have that
\begin{equation}
\nabla_{G}\text{\textrm{-deg}}(\mathscr A,B(V))=\prod_{\xi\in\sigma
_{-}(\mathscr A)}\prod_{i=0}^{r}\left(  {\mbox{\rm deg}_{\mathcal{V}_{i}}%
}\right)  ^{m_{i}(\xi)}. \label{eq:lin-GdegGrad}%
\end{equation}

\noindent\textbf{Acknowledgement.} The authors are grateful to H-P. Wu for the
GAP programing of the topological invariants. C. Garc\'{\i}a was partially
supported by PAPIIT-UNAM through grant IA105217. I. Berezovik and W.
Krawcewicz acknowledge partial support from National Science Foundation
through grant DMS-1413223. W. Krawcewicz was also supported by Guangzhou
University during his visit it the summer 2017.

\end{document}